**Some relations involving the higher derivatives of the Riemann zeta function**

Donal F. Connon

3 June 2015

## Abstract

We show that the higher derivatives of the Riemann zeta function may be expressed in terms of integrals involving the digamma function. Related integrals for the Stieltjes constants are also shown. We also present a formula for $\varsigma^{(n)}(0)$ entirely in terms of the Lehmer constants $b_n$.

## 1. The de Bruijn integral

Using Ramanujan's master theorem, it is shown in Edwards' book [23, p.223] that for $0 < \text{Re}(s) < 1$ the Riemann zeta function may be expressed as the integral

(1.1) $$\varsigma(s) = \frac{\sin(\pi s)}{\pi} \int_0^\infty \frac{\log u - \psi(1+u)}{u^s} du$$

where $\psi(u)$ is the digamma function, which is the logarithmic derivative of the gamma function $\psi(u) = \frac{d}{du} \log \Gamma(u)$. Further information on Ramanujan's master theorem may be found in [2], in Ramanujan's first quarterly report [7, p.298] and in Hardy's lectures on Ramanujan [24, p.186]. Two other derivations of formula (1.1) are given in Titchmarsh's book [30, p.25 & 29], the second one being based on the Müntz formula.

Incidentally, it may be noted that Titchmarsh [30, p.25] has shown how the representation of the Riemann zeta function in (1.1) may also be employed to derive the functional equation for $\varsigma(s)$. As we shall see towards the end of this section, the argument employed by Titchmarsh may be reversed to give another derivation of the integral (1.1).

We may write the above integral as

$$\int_0^\infty \frac{\log u - \psi(1+u)}{u^s} du = \int_0^\infty \frac{\log(1+u) - \psi(1+u)}{u^s} du - \int_0^\infty \frac{\log(1+u) - \log u}{u^s} du$$

and, using integration by parts, we have

$$\int_0^\infty \frac{\log(1+u) - \log u}{u^s} du = \frac{1}{1-s} \log\left(1 + \frac{1}{u}\right) u^{1-s} \Big|_0^\infty + \frac{1}{1-s} \int_0^\infty \frac{u^{-s}}{1+u} du$$

$$= \frac{1}{1-s} \int_0^\infty \frac{u^{-s}}{1+u} du$$

We have the well-known integral [3, p.10]

$$\int_0^\infty \frac{u^{p-1}}{1+u} du = \frac{\pi}{\sin(\pi p)}$$

which is valid for 0 < $p$ < 1. Letting $p = 1-s$ so that for 0 < $s$ < 1 we obtain

(1.2) $$\int_0^\infty \frac{u^{-s}}{1+u} du = \frac{\pi}{\sin(\pi - \pi s)} = \frac{\pi}{\sin(\pi s)}$$

and hence we have determined that for $0 <$ $s$ $< 1$

$$\int_0^\infty \frac{\log(1+u) - \log u}{u^s} du = \frac{\pi}{(1-s)\sin(\pi s)}$$

Therefore we have

(1.3) $$\varsigma(s) = \frac{1}{s-1} + \frac{\sin(\pi s)}{\pi} \int_0^\infty \frac{\log(1+u) - \psi(1+u)}{u^s} du$$

This integral was used by de Bruijn [21] in 1937 to derive (1.1) and it is also valid for 0 < $\text{Re}(s)$ < 1; this paper was de Bruijn's first publication at the age of 18 in response to a problem set by Kloosterman [22] (who had previously derived (1.1) in 1922).This integral is reported in [28, p.102].

Reverting back to (1.1), since $\varsigma(0) = -\frac{1}{2}$ , we see that

(1.4) $$\varsigma(s) - \varsigma(0) = \frac{\sin(\pi s)}{\pi} \int_0^\infty \left[ \log u - \psi(1+u) + \frac{1}{2(1+u)} \right] \frac{du}{u^s}$$

where we have employed (1.2). We write this as

$$\frac{\varsigma(s) - \varsigma(0)}{s} = \frac{\sin(\pi s)}{\pi s} \int_0^\infty \left[ \log u - \psi(1+u) + \frac{1}{2(1+u)} \right] \frac{du}{u^s}$$

and, using L'Hôpital's rule, in the limit as $s \to 0$ we obtain

(1.5) $$\varsigma'(0) = \int_0^\infty \left[ \log u - \psi(1+u) + \frac{1}{2(1+u)} \right] du$$

or equivalently

(1.6) $$\frac{1}{2}\log(2\pi) = \int_0^\infty \left[ \psi(1+u) - \log u - \frac{1}{2(1+u)} \right] du$$

This integral was previously obtained by Berndt and Dixit [8] in a different manner in 2009. Other derivations of this are contained in [18].

More generally, using Leibniz's rule, differentiation of (1.4) results in

(1.7)
$$\varsigma^{(n)}(s) = \frac{1}{\pi}\sum_{k=0}^{n} \binom{n}{k} \pi^k \sin\left(\pi s + \frac{k\pi}{2}\right) \int_0^\infty \left[ \log u - \psi(1+u) + \frac{1}{2(1+u)} \right] \frac{(-1)^{n-k}\log^{n-k} u}{u^s} du$$

where we have noted that

$$\frac{d}{ds}\sin(\pi s) = \pi\cos(\pi s) = \pi\sin\left(\pi s + \frac{\pi}{2}\right)$$

and thus

$$\frac{d^k}{ds^k}\sin(\pi s) = \pi^k \sin\left(\pi s + \frac{k\pi}{2}\right)$$

so that

(1.8) $$\varsigma^{(n)}(0) = \frac{1}{\pi}\sum_{k=0}^{n} \binom{n}{k} \pi^k \sin\left(\frac{k\pi}{2}\right) \int_0^\infty \left[ \log u - \psi(1+u) + \frac{1}{2(1+u)} \right] (-1)^{n-k}\log^{n-k} u\, du$$

For example, we have for $n = 2$

(1.9) $$\varsigma''(0) = 2\int_0^\infty \left[ \psi(1+u) - \log u - \frac{1}{2(1+u)} \right] \log u\, du$$

Ramanujan [6] showed that

(1.10) $$\varsigma''(0) = \gamma_1 + \frac{1}{2}\gamma^2 - \frac{1}{4}\varsigma(2) - \frac{1}{2}\log^2(2\pi)$$

where $\gamma_n$ are the Stieltjes constants [26, p.4]. The higher derivatives $\varsigma^{(n)}(0)$ are also addressed in Apostol's paper [4] where it is shown that they may be expressed in terms involving the Stieltjes constants. Another derivation of (1.9) is shown below in section 3.

The Stieltjes constants $\gamma_n$ are the coefficients of the Laurent expansion of the Riemann zeta function $\varsigma(s)$ about $s = 1$

(1.11) $$\varsigma(s)=\frac{1}{s-1}+\sum_{n=0}^{\infty}\frac{(-1)^n}{n!}\gamma_n(s-1)^n$$

Since $\lim_{s\to 1}\left[\varsigma(s)-\frac{1}{s-1}\right]=\gamma$ it is clear that $\gamma_0=\gamma$. It may be shown, as in [26, p.4], that

(1.12) $$\gamma_n=\lim_{N\to\infty}\left[\sum_{k=1}^{N}\frac{\log^n k}{k}-\frac{\log^{n+1}N}{n+1}\right]=\lim_{N\to\infty}\left[\sum_{k=1}^{N}\frac{\log^n k}{k}-\int_1^N\frac{\log^n t}{t}dt\right]$$

It is easily seen from the Laurent expansion that

(1.12.1) $$\left.\frac{d^n}{ds^n}(s-1)\varsigma(s)\right|_{s=1}=(-1)^{n-1}n\gamma_{n-1}$$

Another example of (1.8) is set out below

$$\varsigma^{(3)}(0)=3\int_0^{\infty}\left[\log u-\psi(1+u)+\frac{1}{2(1+u)}\right]\log^2 u\,du-\pi^2\int_0^{\infty}\left[\log u-\psi(1+u)+\frac{1}{2(1+u)}\right]du$$

and substituting (1.6) we have

(1.13) $$\varsigma^{(3)}(0)=3\int_0^{\infty}\left[\log u-\psi(1+u)+\frac{1}{2(1+u)}\right]\log^2 u\,du+\frac{\pi^2}{2}\log(2\pi)$$

□

Using Euler's reflection formula

$$\Gamma(s)\Gamma(1-s)=\frac{\pi}{\sin\pi s}$$

we may also express (1.4) as

$$\varsigma(s)-\varsigma(0)=\frac{1}{\Gamma(s)\Gamma(1-s)}\int_0^{\infty}\left[\log u-\psi(1+u)+\frac{1}{2(1+u)}\right]\frac{du}{u^s}$$

and, for convenience, we denote $F(s,u)$ as

$$F(s,u)=\frac{1}{\Gamma(s)\Gamma(1-s)u^s}$$

We need to deal with the troublesome factor of $\Gamma(s)$ in the denominator and, to this end, we write this in the equivalent form

$$F(s,u) = \frac{su^{-s}}{\Gamma(1+s)\Gamma(1-s)}$$

First of all, we employ the Leibniz differentiation formula to obtain

$$\frac{\partial^n}{\partial s^n} F(s,u) = \sum_{k=0}^{n} \binom{n}{k} \frac{\partial^{n-k}}{\partial s^{n-k}} [su^{-s}] \frac{d^k}{ds^k} \frac{1}{\Gamma(1+s)\Gamma(1-s)}$$

We see that

$$\frac{d}{ds} \frac{1}{\Gamma(1+s)\Gamma(1-s)} = \frac{1}{\Gamma(1+s)\Gamma(1-s)} [\psi(1-s) - \psi(1+s)]$$

and hence we have

(1.14) $$\frac{d^k}{ds^k} \frac{1}{\Gamma(1+s)\Gamma(1-s)}$$

$$= \frac{1}{\Gamma(1+s)\Gamma(1-s)} Y_k \left(\psi(1-s) - \psi(1+s), -\psi^{(1)}(1-s) - \psi^{(1)}(1+s), ..., (-1)^{k-1}\psi^{(k-1)}(1-s) - \psi^{(k-1)}(1+s)\right)$$

in terms of the (exponential) complete Bell polynomials $Y_n(x_1,...,x_n)$ which are defined by $Y_0 = 1$ and for $n \ge 1$

(1.15) $$Y_n(x_1,...,x_n) = \sum_{\pi(n)} \frac{n!}{k_1!\, k_2! ... k_n!} \left(\frac{x_1}{1!}\right)^{k_1} \left(\frac{x_2}{2!}\right)^{k_2} ... \left(\frac{x_n}{n!}\right)^{k_n}$$

where the sum is taken over all partitions $\pi(n)$ of $n$, i.e. over all sets of integers $k_j$ such that

$$k_1 + 2k_2 + 3k_3 + ... + nk_n = n$$

The complete Bell polynomials have integer coefficients and the first six are set out below (Comtet [15, p.307])

$$Y_1(x_1) = x_1$$

$$Y_2(x_1,x_2) = x_1^2 + x_2$$

$$Y_3(x_1,x_2,x_3) = x_1^3 + 3x_1x_2 + x_3$$

$$Y_4(x_1,x_2,x_3,x_4) = x_1^4 + 6x_1^2x_2 + 4x_1x_3 + 3x_2^2 + x_4$$

$$Y_5(x_1,x_2,x_3,x_4,x_5) = x_1^5 + 10x_1^3x_2 + 10x_1^2x_3 + 15x_1x_2^2 + 5x_1x_4 + 10x_2x_3 + x_5$$

$$Y_6(x_1,x_2,x_3,x_4,x_5,x_6) = x_1^6 + 6x_1x_5 + 15x_2x_4 + 10x_2^3 + 15x_1^2x_4 + 15x_2^3 + 60x_1x_2x_3$$

$$+20x_1^3x_3 + 45x_1^2x_2^2 + 15x_1^4x_1 + x_6$$

All of the information that we require for the purposes of this paper regarding the Bell polynomials is contained in [17]; more detailed expositions may be found in [10] and [15].

In particular, suppose that $h'(x) = h(x)g(x)$ then we have

$$h^{(r)}(x) = h(x)Y_r\left(g(x), g^{(1)}(x), ..., g^{(r-1)}(x)\right)$$

which is used in (1.14) above and elsewhere in this paper.

We also have

$$\frac{\partial^j}{\partial s^j}[su^{-s}] = (-1)^j\left[su^{-s}\log^j u - ju^{-s}\log^{j-1} u\right]$$

so that

$$\frac{\partial^n}{\partial s^n}F(s,v) = \frac{u^{-s}}{\Gamma(1+s)\Gamma(1-s)}\sum_{k=0}^{n}\binom{n}{k}(-1)^{n-k}\left(s\log^{n-k}u - (n-k)\log^{n-k-1}u\right)\mathbf{Y}_k(\boldsymbol{\psi}(s))$$

where, for convenience, we denote

$$\mathbf{Y}_k(\boldsymbol{\psi}(s)) = Y_k\left(\psi(1-s) - \psi(1+s), -\psi^{(1)}(1-s) - \psi^{(1)}(1+s), ..., (-1)^{k-1}\psi^{(k-1)}(1-s) - \psi^{(k-1)}(1+s)\right)$$

When $s = 0$ this becomes

$$\frac{\partial^n}{\partial s^n}F(0,v) = \sum_{k=0}^{n}\binom{n}{k}(-1)^{n-k+1}(n-k)\log^{n-k-1}v\cdot\mathbf{Y}_k(\boldsymbol{\psi}(0))$$

$$= \sum_{k=0}^{n-1}\binom{n}{k}(-1)^{n-k+1}(n-k)\log^{n-k-1}u\cdot\mathbf{Y}_k(\boldsymbol{\psi}(0))$$

where

$$\mathbf{Y}_k(\boldsymbol{\psi}(0)) = Y_k\left(0, -2\psi^{(1)}(1), ..., [(-1)^{k-1} - 1]\psi^{(k-1)}(1)\right)$$

We note that [28, p.22]

$$\psi^{(r)}(a) = (-1)^{r+1}r!\varsigma(r+1,a)$$

and hence we see that

$$\mathbf{Y}_k(\boldsymbol{\psi}(0)) = Y_k\left(0, -2\varsigma(2)1!, ..., [(-1)^{k+1} - 1]\varsigma(k)(k-1)!\right)$$

Using the elementary binomial identity

$$(n-k)\binom{n}{k} = n\binom{n-1}{k}$$

we obtain

$$\frac{\partial^n}{\partial s^n} F(0,v) = n\sum_{k=0}^{n-1}\binom{n-1}{k}(-1)^{n-k+1}\log^{n-k-1} u \cdot \mathbf{Y}_k(\boldsymbol{\psi}(0))$$

We therefore conclude that

(1.16)

$$\varsigma^{(n)}(0) = n(-1)^n \sum_{k=0}^{n-1}\binom{n-1}{k}(-1)^{k+1}\mathbf{Y}_k(\boldsymbol{\psi}(0))\int_0^\infty \left[\log u - \psi(1+u) + \frac{1}{2(1+u)}\right]\log^{n-k-1} u\, du$$

With $n = 2$ we see that

$$\varsigma''(0) = -2\int_0^\infty \left[\log u - \psi(1+u) + \frac{1}{2(1+u)}\right]\log u\, du$$

$$+\mathbf{Y}_1(\boldsymbol{\psi}(0))\int_0^\infty \left[\log u - \psi(1+u) + \frac{1}{2(1+u)}\right] du$$

and this concurs with (1.9) since $\mathbf{Y}_1(\boldsymbol{\psi}(0)) = 0$.

We note that the factor of $n$ in (1.16) suggests that $\varsigma^{(n)}(0)$ is unbounded as $n \to \infty$ which is consistent with the known limit $\lim_{n\to\infty} \frac{\varsigma^{(n)}(0)}{n!} = -1$ mentioned by Apostol [4].

□

We may write (1.4) as

$$[\varsigma(s) - \varsigma(0)]h(s) = s\int_0^\infty \left[\log u - \psi(1+u) + \frac{1}{2(1+u)}\right]\frac{du}{u^s}$$

where $h(s) = \Gamma(1+s)\Gamma(1-s)$. We see that

$$h'(s) = h(s)[\psi(1+s) - \psi(1-s)]$$

and the Leibniz rule gives us for $s = 0$

(1.16.1) $$\sum_{k=0}^{n-1}\binom{n}{k}\varsigma^{(n-k)}(0)\mathbf{Y}_k(-\boldsymbol{\psi}(0)) = n(-1)^{n-1}\int_0^\infty\left[\log u - \psi(1+u) + \frac{1}{2(1+u)}\right]\log^{n-1} u\, du$$

where $\mathbf{Y}_k(-\boldsymbol{\psi}(0)) = Y_k\left(0, 2\varsigma(2)1!, ..., -[(-1)^{k+1}-1]\varsigma(k)(k-1)!\right)$.

□

In this section we reverse the argument employed by Titchmarsh [28, p.25] to give another derivation of (1.1).

We see from (1.2) that for $0 < p < 1$

$$\int_0^\infty \frac{t^{-p}}{1+t}dt = \frac{\pi}{\sin(\pi p)}$$

and the substitution $t = v/x$ gives us

(1.17) $$\int_0^\infty \frac{v^{-p}}{x+v}dv = \frac{\pi}{x^p\sin(\pi p)}$$

Letting $v = u^2$ results in

$$\int_0^\infty \frac{u^{1-2p}}{x+u^2}du = \frac{\pi}{2x^p\sin(\pi p)}$$

and with $x \to x^2$ we see that

$$\int_0^\infty \frac{u^{1-2p}}{x^2+u^2}du = \frac{\pi}{2x^{2p}\sin(\pi p)}$$

We now designate $s = 2p-1$ and write this for $-1 < s < 1$ as

(1.18) $$\int_0^\infty \frac{u^{-s}}{x^2+u^2}du = \frac{\pi}{2x^{s+1}\cos(\pi s/2)}$$

We have the well known integral for the digamma function [32, p.251], another proof of which is shown below in (3.8)

$$\psi(u) = -\frac{1}{2u} + \log u - 2\int_0^\infty \frac{x}{(x^2+u^2)(e^{2\pi x}-1)}dx$$

which we may write as

$$\log u-\psi(1+u)=-\frac{1}{2u}+2\int_0^\infty \frac{x}{(x^2+u^2)(e^{2\pi x}-1)}dx$$

Since

(1.19) $$\int_0^\infty \frac{dx}{x^2+u^2}=\frac{1}{u}\tan^{-1}(x/u)\Big|_0^\infty=\frac{\pi}{2u}$$

we have

(1.20) $$\log u-\psi(1+u)=2\int_0^\infty \frac{x}{x^2+u^2}\left[\frac{1}{e^{2\pi x}-1}-\frac{1}{2\pi x}\right]dx$$

It may be noted that (1.19) may be deduced by letting $s=0$ in (1.18).

We now multiply (1.20) by $u^{-s}$ and integrate over $[0,\infty)$ to obtain

$$\int_0^\infty \frac{\log u-\psi(1+u)}{u^s}du=2\int_0^\infty\int_0^\infty \frac{xu^{-s}}{x^2+u^2}\left[\frac{1}{e^{2\pi x}-1}-\frac{1}{2\pi x}\right]dx\,du$$

$$=2\int_0^\infty \left[\frac{1}{e^{2\pi x}-1}-\frac{1}{2\pi x}\right]x\,dx\int_0^\infty \frac{u^{-s}}{x^2+u^2}du$$

$$=\frac{\pi}{\cos(\pi s/2)}\int_0^\infty \left[\frac{1}{e^{2\pi x}-1}-\frac{1}{2\pi x}\right]x^{(1-s)-1}\,dx$$

where we have used (1.18) in the final part.

We note that

$$\int_0^\infty \left[\frac{1}{e^{2\pi x}-1}-\frac{1}{2\pi x}\right]x^{(1-s)-1}\,dx=\frac{1}{(2\pi)^{1-s}}\int_0^\infty \left[\frac{1}{e^t-1}-\frac{1}{t}\right]t^{(1-s)-1}\,dt$$

It is well known that for $\mathrm{Re}(s)>1$ [28, p.96]

$$\varsigma(s)\Gamma(s)=\int_0^\infty \frac{t^{s-1}}{e^t-1}dt$$

and we have for $0<\mathrm{Re}(s)<1$ ([28, p.23] and [29, p.162])

$$\varsigma(s)\Gamma(s)=\int_0^\infty \left[\frac{1}{e^t-1}-\frac{1}{t}\right]t^{s-1}dt$$

This then gives us

$$\varsigma(1-s)\Gamma(1-s)=\int_0^\infty\left[\frac{1}{e^t-1}-\frac{1}{t}\right]t^{-s}dt$$

and we see that

$$\int_0^\infty \frac{\log u-\psi(1+u)}{u^s}du=\frac{\pi}{\cos(\pi s/2)}(2\pi)^{s-1}\varsigma(1-s)\Gamma(1-s)$$

With the functional equation for the Riemann zeta function

(1.21) $\varsigma(1-s)=2(2\pi)^{-s}\Gamma(s)\cos(\pi s/2)\varsigma(s)$

and employing Euler's reflection formula for the gamma function we obtain (1.1)

$$\int_0^\infty \frac{\log u-\psi(1+u)}{u^s}du=\frac{\pi}{\sin\pi s}\varsigma(s)$$

□

We now multiply (1.20) by $u^{-s}\log u$ and integrate over $[0,\infty)$ to obtain

$$\int_0^\infty \frac{\log u-\psi(1+u)}{u^s}\log u\,du=2\int_0^\infty\int_0^\infty \frac{xu^{-s}\log u}{x^2+u^2}\left[\frac{1}{e^{2\pi x}-1}-\frac{1}{2\pi x}\right]dx\,du$$

$$=2\int_0^\infty\left[\frac{1}{e^{2\pi x}-1}-\frac{1}{2\pi x}\right]x\,dx\int_0^\infty \frac{u^{-s}\log u}{x^2+u^2}du$$

We recall (1.18)

$$\int_0^\infty \frac{u^{-s}}{x^2+u^2}du=\frac{\pi}{2x^{s+1}\cos(\pi s/2)}$$

Carslaw [9, p.212] indicates that differentiation under the integral sign is valid here and we obtain

(1.22) $$\int_0^\infty \frac{u^{-s}\log u}{x^2+u^2}du=\frac{\pi}{2}x^{-(s+1)}\frac{(\pi/2)\sin(\pi s/2)-\cos(\pi s/2)\log x}{\cos^2(\pi s/2)}$$

Hence we obtain

$$\int_0^\infty \frac{\log u-\psi(1+u)}{u^s}\log u\,du$$

$$=\frac{\pi}{\cos^2(\pi s/2)}\int_0^\infty\left[\frac{1}{e^{2\pi x}-1}-\frac{1}{2\pi x}\right]\frac{(\pi/2)\sin(\pi s/2)-\cos(\pi s/2)\log x}{x^s}dx$$

which may be equated with the result obtained by differentiating (1.1).

We note that letting $s=0$ in (1.22) gives us

(1.23) $$\int_0^\infty \frac{\log u}{x^2+u^2}du = \frac{\pi}{2x}\log x$$

and this is a particular case of (3.7) below.

**2. An integral representation of the Stieltjes constants**

With reference to (1.3), which is valid for $0<\operatorname{Re}(s)<1$, we note that integration by parts gives us

$$\int_0^\infty \frac{\log(1+u)-\psi(1+u)}{u^s}du = [\log(1+u)-\psi(1+u)]\frac{u^{1-s}}{1-s}\bigg|_0^\infty - \frac{1}{1-s}\int_0^\infty \left[\frac{1}{1+u}-\psi'(1+u)\right]\frac{1}{u^{s-1}}du$$

and, since $\lim_{u\to\infty}[\log(1+u)-\psi(1+u)]=0$, we find that the integrated parts vanish and thus

$$\int_0^\infty \frac{\log(1+u)-\psi(1+u)}{u^s}du = -\frac{1}{1-s}\int_0^\infty \left[\frac{1}{1+u}-\psi'(1+u)\right]\frac{1}{u^{s-1}}du$$

Hence we obtain de Bruijn's formula [21] which is valid for $0<\operatorname{Re}(s)<2; s\neq 1$

(2.1) $$\varsigma(s) = \frac{1}{s-1} - \frac{\sin(\pi s)}{\pi(s-1)}\int_0^\infty \left[\psi'(1+u)-\frac{1}{1+u}\right]\frac{1}{u^{s-1}}du$$

and this corrects two misprints in [28, p.103]. We note that (2.1) continues to be recorded incorrectly in Choi's recent paper [11].

Since $\lim_{s\to 1}\frac{\sin(\pi s)}{\pi(s-1)} = \lim_{s\to 1}\frac{\pi\cos(\pi s)}{\pi} = -1$ we see from (2.1) that

$$\lim_{s\to 1}\left[\varsigma(s)-\frac{1}{s-1}\right] = \int_0^\infty \left[\psi'(1+u)-\frac{1}{1+u}\right]du$$

$$= [\psi(1+u)-\log(1+u)]\Big|_0^\infty$$

Since $\lim_{u\to\infty}[\psi(1+u)-\log(1+u)]=0$ we obtain

(2.2) $$\gamma = \int_0^\infty \left[\psi'(1+u)-\frac{1}{1+u}\right]du$$

This is incorrectly reported in [11]. The integral is in accordance with the well known limit

$$\gamma = \lim_{s\to 1}\left[\varsigma(s) - \frac{1}{s-1}\right]$$

We write (2.1) as

$$\varsigma(s) - \frac{1}{s-1} = -\frac{\sin(\pi s)}{\pi(s-1)}\int_0^\infty \left[\psi'(1+u) - \frac{1}{1+u}\right]\frac{1}{u^{s-1}}\,du$$

and then employ $\frac{d^n}{ds^n}\left[\varsigma(s) - \frac{1}{s-1}\right]\Bigg|_{s=1} = (-1)^n \gamma_n$ to evaluate the Stieltjes constants (which is the approach adopted by Choi in [11]). However, prior to the publication of Choi's paper [11] in 2013 I carried out the analysis in a slightly different manner.

In this alternative approach, we write (2.1) as

$$\varsigma(s) - \frac{1}{s-1} = -\frac{1}{(s-1)\Gamma(s)\Gamma(1-s)}\int_0^\infty \left[\psi'(1+u) - \frac{1}{1+u}\right]\frac{1}{u^{s-1}}\,du$$

$$= \frac{1}{\Gamma(s)\Gamma(2-s)}\int_0^\infty \left[\psi'(1+u) - \frac{1}{1+u}\right]\frac{1}{u^{s-1}}\,du$$

With $h(s) = \frac{1}{\Gamma(s)\Gamma(2-s)}$, we find that $h(s) = h(s)[\psi(2-s) - \psi(s)]$ and we therefore obtain

$$\frac{d^n}{ds^n}\left[\varsigma(s) - \frac{1}{s-1}\right] = h(s)\sum_{k=0}^n \binom{n}{k} \mathbf{Y}_k(\boldsymbol{\psi}_1(s))\int_0^\infty \left[\psi'(1+u) - \frac{1}{1+u}\right](-1)^{n-k}\log^{n-k} u\,du$$

where

$$\mathbf{Y}_k(\boldsymbol{\psi}_1(s)) \equiv Y_k\left(\psi(2-s) - \psi(s), -\psi'(2-s) - \psi'(s), ..., (-1)^{k-1}\psi^{(k-1)}(2-s) - \psi^{(k-1)}(s)\right)$$

With $s = 1$ we have

$$(2.3) \qquad \gamma_n = \sum_{k=0}^n \binom{n}{k}(-1)^k\, \mathbf{Y}_k(\boldsymbol{\psi}_1(1))\int_0^\infty \left[\psi'(1+u) - \frac{1}{1+u}\right]\log^{n-k} u\,du$$

where, using $\psi^{(r)}(a) = (-1)^{r+1} r!\varsigma(r+1,a)$, we have

$$\mathbf{Y}_k(\boldsymbol{\psi}_1(1)) = Y_k\left(0, -2\varsigma(2)1!, ..., -[1+(-1)^k]\varsigma(k)(k-1)!\right)$$

For example, with $n = 1$ we obtain

(2.4) $$\gamma_1 = \int_0^\infty \left[ \psi'(1+u) - \frac{1}{1+u} \right] \log u \, du$$

□

We may also obtain an expression for the individual integrals contained in (2.3) by writing (2.1) in the following format

$$\left[ \varsigma(s) - \frac{1}{s-1} \right] \Gamma(s)\Gamma(2-s) = \int_0^\infty \left[ \psi'(1+u) - \frac{1}{1+u} \right] \frac{1}{u^{s-1}} \, du$$

Defining $f(s)$ as $f(s) = \Gamma(s)\Gamma(2-s)$, the Leibniz differentiation formula evaluated at $s=1$ results in

$$(-1)^n \int_0^\infty \left[ \psi'(1+u) - \frac{1}{1+u} \right] \log^n u \, du = \sum_{k=0}^n \binom{n}{k} f^{(k)}(1)(-1)^{n-k} \gamma_{n-k}$$

where we have employed $\left. \frac{d^n}{ds^n} \left[ \varsigma(s) - \frac{1}{s-1} \right] \right|_{s=1} = (-1)^n \gamma_n$.

We note that the derivative of $f(s)$ is

$$f'(s) = \Gamma(s)\Gamma(2-s)[\psi(s) - \psi(2-s)]$$

which results in

$$f^{(k)}(s) = f(s) Y_k \left( \psi(s) - \psi(2-s), \psi'(s) + \psi'(2-s), ..., \psi^{(k-1)}(s) - (-1)^{k-1} \psi^{(k-1)}(2-s) \right)$$

and hence we obtain

$$f^{(k)}(1) = Y_k \left( 0, 2\psi'(1), ..., [1+(-1)^k] \psi^{(k-1)}(1) \right)$$

We note that [28, p.22]

$$\psi^{(r)}(a) = (-1)^{r+1} r! \varsigma(r+1, a)$$

which gives us

$$f^{(k)}(1) = Y_k \left( 0, 2.1! \varsigma(2), ..., [1+(-1)^k](k-1)! \varsigma(k) \right)$$

Hence we obtain the integral

(2.5) $$\int_0^\infty \left[ \psi'(1+u) - \frac{1}{1+u} \right] \log^n u \, du = \sum_{k=0}^n \binom{n}{k} \mathbf{Y}_k (-1)^k \gamma_{n-k}$$

where $\mathbf{Y}_k \equiv Y_k\left(0, 2.1!\varsigma(2), ..., [1+(-1)^k](k-1)!\varsigma(k)\right)$. It may be noted that

$$f^{(k)}(1) = \mathbf{Y}_k(-\boldsymbol{\psi}_1(1)) = Y_k\left(0, 2\varsigma(2)1!, ..., [(-1)^k+1]\varsigma(k)(k-1)!\right)$$

□

Integration by parts gives us

$$\int \left[\psi'(1+u) - \frac{1}{1+u}\right] \log u\, du = \left[\psi(1+u) - \log(1+u)\right] \log u - \int \frac{\psi(1+u) - \log(1+u)}{u}\, du$$

and this gives us the definite integral

$$\int_1^\infty \left[\psi'(1+u) - \frac{1}{1+u}\right] \log u\, du = -\int_1^\infty \frac{\psi(1+u) - \log(1+u)}{u}\, du$$

We write the latter integral as

$$\int_1^\infty \frac{\psi(1+u) - \log(1+u)}{u}\, du = \int_1^\infty \frac{\psi(1+u) - \log u}{u}\, du + \int_1^\infty \frac{\log u - \log(1+u)}{u}\, du$$

We note from [28, p.106] that

$$\frac{\log u - \log(1+u)}{u} = -\frac{d}{du} Li_2\left(-\frac{1}{u}\right)$$

and integration results in

$$\int_1^x \frac{\log u - \log(1+u)}{u}\, du = Li_2(-1) - Li_2\left(-\frac{1}{x}\right)$$

in terms of the polylogarithm function.

Hence we have

$$(2.6) \qquad \int_1^\infty \frac{\log u - \log(1+u)}{u}\, du = -\frac{1}{2}\varsigma(2)$$

and thus we obtain

$$\int_1^\infty \left[\psi'(1+u) - \frac{1}{1+u}\right] \log u\, du = -\int_1^\infty \frac{\psi(1+u) - \log u}{u}\, du + \frac{1}{2}\varsigma(2)$$

We showed in Eq.(3.35.1) in [19] that

(2.7) $$\int_0^1 \frac{\psi(1+u)+\gamma}{u}\,du+\int_1^\infty \frac{\psi\left(1+u\right)-\log u}{u}\,du=\varsigma(2)-\gamma_1$$

and it is known [24 , p.142] that

(2.8) $$\int_0^1 \frac{\psi(1+u)+\gamma}{u}\,du=\sum_{n=1}^\infty \frac{\log(n+1)}{n(n+1)}$$

It is not known whether the above series has a representation in closed form.

Hence we have the integrals

(2.9) $$\int_1^\infty \left[\psi'(1+u)-\frac{1}{1+u}\right]\log u\,du=\sum_{n=1}^\infty \frac{\log(n+1)}{n(n+1)}+\gamma_1-\frac{1}{2}\varsigma(2)$$

and

(2.10) $$\int_0^1 \left[\psi'(1+u)-\frac{1}{1+u}\right]\log u\,du=\frac{1}{2}\varsigma(2)-\sum_{n=1}^\infty \frac{\log(n+1)}{n(n+1)}$$

Since

$$\int \frac{\log u}{1+u}\,du=Li_2(-u)+\log u\log(1+u)$$

we have

$$\int_0^1 \frac{\log u}{1+u}\,du=-\frac{1}{2}\varsigma(2)$$

and we deduce that

(2.11) $$\int_0^1 \psi'(1+u)\log u\,du=-\sum_{n=1}^\infty \frac{\log(n+1)}{n(n+1)}$$

This may also be deduced by integrating (2.8) by parts. Indeed, it is easy to directly evaluate the integral be by noting that

$$\psi'(1+u)=\varsigma(2,1+u)=\sum_{n=0}^\infty \frac{1}{(n+u+1)^2}$$

and using

$$\int \frac{\log t}{(t+a)^2}\,dt=\frac{t\log t-(t+a)\log(t+a)}{a(t+a)}$$

whereupon we obtain the equivalent version

(2.12) $$\int_0^1 \psi'(1+u)\log u\, du = -\sum_{n=1}^{\infty}\frac{1}{n}\log\frac{n+1}{n}$$

It may be noted that Cohen [24, p.142] has also stated that

$$\sum_{n=1}^{\infty}\frac{\log(n+1)}{n(n+1)} = \int_0^1 \frac{(1-x)\log(1-x)}{x\log x}dx$$

$$= \sum_{n=1}^{\infty}(-1)^{n+1}\frac{\varsigma(n+1)}{n}$$

$$= -\sum_{n=2}^{\infty}\varsigma'(n)$$

$$= \sum_{n=1}^{\infty}\frac{1}{n}\log\frac{n+1}{n}$$ □

We would emphasise that it is not necessary to employ the complete Bell polynomials; for example, using the Leibniz rule to differentiate

$$(s-1)\varsigma(s) = 1 - \frac{\sin(\pi s)}{\pi}\int_0^{\infty}\left[\psi'(1+u) - \frac{1}{1+u}\right]\frac{1}{u^{s-1}}du$$

we obtain

$$\frac{d^n}{ds^n}\left[(s-1)\varsigma(s)\right]\bigg|_{s=1} = -\sum_{j=0}^{n}\binom{n}{j}\pi^{n-j-1}\sin\left(\frac{(n-j)\pi}{2}\right)\int_0^{\infty}\left[\psi'(1+u) - \frac{1}{1+u}\right](-1)^j\log^j u\, du$$

where we have used

$$\frac{d^k}{ds^k}\sin(\pi s) = \pi^k\sin\left(\pi s + \frac{k\pi}{2}\right)$$

$$\frac{d^k}{ds^k}\sin(\pi s)\bigg|_{s=1} = \pi^k\sin\left(\frac{k\pi}{2}\right)$$

Hence we obtain

(2.13) $$(-1)^n n\gamma_{n-1} = \sum_{j=0}^{n}\binom{n}{j}\pi^{n-j-1}\sin\left(\frac{(n-j)\pi}{2}\right)\int_0^{\infty}\left[\psi'(1+u) - \frac{1}{1+u}\right](-1)^j\log^j u\, du$$

## 3. The Hurwitz zeta function

Adamchik [1] noted that the Hermite integral for the Hurwitz zeta function may be derived from the Abel-Plana summation formula

(3.1) $$\sum_{k=0}^{\infty} f(k) = \frac{1}{2} f(0) + \int_0^{\infty} f(x)\,dx + i\int_0^{\infty} \frac{f(ix) - f(-ix)}{e^{2\pi x} - 1}\,dx$$

which applies to functions which are analytic in the right-hand plane and satisfy the convergence condition $\lim_{y\to\infty} e^{-2\pi y} \left| f(x+iy) \right| = 0$ uniformly on any finite interval of $x$. Derivations of the Abel-Plana summation formula may be found in [5], [31, p.108] and [32, p.145].

Letting $f(k) = (k+u)^{-s}$ we obtain

(3.2) $$\varsigma(s,u) = \sum_{k=0}^{\infty} \frac{1}{(k+u)^s} = \frac{u^{-s}}{2} + \frac{u^{1-s}}{s-1} + i\int_0^{\infty} \frac{(u+ix)^{-s} - (u-ix)^{-s}}{e^{2\pi x} - 1}\,dx$$

Then, noting that

$$(u+ix)^{-s} - (u-ix)^{-s} = (re^{i\theta})^{-s} - (re^{-i\theta})^{-s}$$

$$= r^{-s}[e^{-is\theta} - e^{is\theta}]$$

$$= \frac{2}{i(u^2+x^2)^{s/2}} \sin(s\tan^{-1}(x/u))$$

we may write (3.2) as Hermite's integral for $\varsigma(s,u)$

(3.3) $$\varsigma(s,u) = \frac{u^{-s}}{2} + \frac{u^{1-s}}{s-1} + 2\int_0^{\infty} \frac{\sin(s\tan^{-1}(x/u))}{(u^2+x^2)^{s/2}(e^{2\pi x} - 1)}\,dx$$

We now take one step back and differentiate the intermediate equation (3.2) with respect to $s$ to obtain

(3.4)
$$\varsigma'(s,u) = -\frac{1}{2} u^{-s}\log u - \frac{u^{1-s}[1+(s-1)\log u]}{(s-1)^2} + i\int_0^{\infty} \frac{(u-ix)^{-s}\log(u-ix) - (u+ix)^{-s}\log(u+ix)}{e^{2\pi x} - 1}\,dx$$

With $s = 0$ in (3.4) we obtain

$$\varsigma'(0,u) = \left(u - \frac{1}{2}\right)\log u - u + i\int_0^{\infty} \frac{\log(u-ix) - \log(u+ix)}{e^{2\pi x} - 1}\,dx$$

This may be written as

$$\varsigma'(0,u) = \left(u - \frac{1}{2}\right)\log u - u + 2\int_0^{\infty} \frac{\tan^{-1}(x/u)}{e^{2\pi x} - 1}\,dx$$

and using Lerch's identity [28, p.92]

$$\log\Gamma(u) = \varsigma'(0,u) + \frac{1}{2}\log(2\pi)$$

we see that this is equivalent to Binet's second formula for $\log\Gamma(u)$ (which is derived in a different manner in, for example, [32, p.251])

$$\text{(3.5)} \qquad \log\Gamma(u) = \left(u - \frac{1}{2}\right)\log u - u + \frac{1}{2}\log(2\pi) + 2\int_0^\infty \frac{\tan^{-1}(x/u)}{e^{2\pi x}-1}dx$$

This formula was also derived by Ramanujan [7, Part II, p.221] in the case where $u$ is a positive integer.

We now consider the second derivative of the Hurwitz zeta function (3.2)

$$\varsigma''(s,u) = \frac{1}{2}u^{-s}\log^2 u + \frac{u^{1-s}[(s-1)\log u + 2]}{(s-1)^3} + \frac{u^{1-s}\log^2 u}{s-1} + \frac{u^{1-s}\log u}{(s-1)^2}$$

$$-i\int_0^\infty \frac{(u-ix)^{-s}\log^2(u-ix) - (u+ix)^{-s}\log^2(u+ix)}{e^{2\pi x}-1}dx$$

where with $s = 0$ we have

$$\text{(3.6)} \qquad \varsigma''(0,u) = \left(\frac{1}{2} - u\right)\log^2 u + 2u\log u - 2u - 2\int_0^\infty \frac{\log\left(u^2+x^2\right)\tan^{-1}(x/u)}{e^{2\pi x}-1}dx$$

Using contour integration, Holland [25, p.191] showed that for $x > 0$

$$\text{(3.7)} \qquad \int_0^\infty \frac{\log t}{(t+u)^2 + x^2}dt = \frac{1}{2x}\log\left(u^2+x^2\right)\tan^{-1}(x/u)$$

This integral is also valid in the limit as $x \to 0$ because L'Hôpital's rule shows that $\lim_{x\to 0}\frac{1}{2x}\log\left(u^2+x^2\right)\tan^{-1}(x/u) = \log u$ which concurs with the integral $\int_0^\infty \frac{\log t}{(t+u)^2}dt = \log u$.

We multiply (3.7) by $\frac{x}{e^{2\pi x}-1}$ and integrate with respect to $x$ to obtain

$$\int_0^\infty \frac{\log\left(u^2+x^2\right)\tan^{-1}(x/u)}{e^{2\pi x}-1}dx = 2\int_0^\infty \frac{x}{e^{2\pi x}-1}dx\int_0^\infty \frac{\log t}{(t+u)^2+x^2}dt$$

$$= 2\int_0^\infty \int_0^\infty \frac{x \log t}{[(t+u)^2 + x^2](e^{2\pi x} - 1)} dx\, dt$$

Differentiating (3.5) results in [28, p.16]

(3.8) $$\psi(u) = -\frac{1}{2u} + \log u - 2\int_0^\infty \frac{x}{(u^2 + x^2)(e^{2\pi x} - 1)} dx$$

and we see that with $u \to t + u$

(3.9) $$\psi(t+u) - \log(t+u) + \frac{1}{2(t+u)} = -2\int_0^\infty \frac{x}{[(t+u)^2 + x^2](e^{2\pi x} - 1)} dx$$

We multiply this by $\log u$ and integrate; this gives us

$$\int_0^\infty \left[ \psi(t+u) - \log(t+u) + \frac{1}{2(t+u)} \right] \log du = -2\int_0^\infty \frac{\log u}{[(t+u)^2 + x^2]} du \int_0^\infty \frac{x}{(e^{2\pi x} - 1)} dx$$

$$= -\int_0^\infty \frac{\log\left(t^2 + x^2\right) \tan^{-1}\left(x/t\right)}{e^{2\pi x} - 1} dx$$

We then obtain

(3.10) $$\int_0^\infty \frac{\log\left(u^2 + x^2\right) \tan^{-1}\left(x/u\right)}{e^{2\pi x} - 1} dx = -\int_0^\infty \left[ \psi(t+u) - \log(t+u) + \frac{1}{2(t+u)} \right] \log t\, dt$$

and using (3.6) we obtain

(3.11)
$$\varsigma''(0,u) = \left(\frac{1}{2} - u\right) \log^2 u + 2u \log u - 2u + 2\int_0^\infty \left[ \psi(t+u) - \log(t+u) + \frac{1}{2(t+u)} \right] \log t\, dt$$

With $u = 1$ we get

(3.12) $$\varsigma''(0) = -2 + 2\int_0^\infty \left[ \psi(t+1) - \log(t+1) + \frac{1}{2(t+1)} \right] \log t\, dt$$

and, as shown below, this may be reconciled with (1.9).

It is easy to determine that

$$\int_0^N \left[ \log u - \log(1+u) + \frac{1}{1+u} \right] \log u\, du = N(\log N - 1)[\log N - \log(1+N)] + \log(1+N)$$

$$= -N(\log N - 1)\log\left(1+\frac{1}{N}\right) + \log(1+N)$$

$$= -N\log N\log\left(1+\frac{1}{N}\right) + N\log\left(1+\frac{1}{N}\right) + \log\left(1+\frac{1}{N}\right) + \log N$$

and, since $\log\left(1+\frac{1}{N}\right) = \frac{1}{N} + O\left(\frac{1}{N^2}\right)$, we then see that

(3.13) $$\int_0^\infty \left[\log u - \log(1+u) + \frac{1}{1+u}\right]\log u\, du = 1$$

We note that

$$\int_0^\infty \left[\psi(1+u) - \log(1+u) + \frac{1}{2(1+u)}\right]\log u\, du$$

$$= \int_0^\infty \left[\psi(1+u) - \log u - \frac{1}{2(1+u)} + \log u - \log(1+u) + \frac{1}{1+u}\right]\log u\, du$$

and thus we have

$$\int_0^\infty \left[\psi(1+u) - \log(1+u) + \frac{1}{2(1+u)}\right]\log u\, du = 1 + \int_0^\infty \left[\psi(1+u) - \log u - \frac{1}{2(1+u)}\right]\log u\, du$$

We therefore see that (1.9) and (3.12) are equivalent.

□

By differentiating (3.11) and noting that [20]

$$\frac{\partial^2}{\partial s^2}\frac{\partial}{\partial u}\varsigma(s,u)\bigg|_{s=0} = 2\gamma_1(u)$$

we obtain

(3.13.1) $$\gamma_1(u) = \frac{1}{2u}\log u - \frac{1}{2}\log^2 u + \int_0^\infty \left[\psi'(t+u) - \frac{1}{t+u} - \frac{1}{2(t+u)^2}\right]\log t\, dt$$

We have

$$\int \frac{\log t}{(t+u)^2}dt = \frac{t\log t - (t+u)\log(t+u)}{u(t+u)}$$

$$= -\frac{1}{u}\log\left(1+\frac{u}{t}\right) - \frac{\log t}{t+u}$$

which gives us the definite integral

$$\int_0^\infty \frac{\log t}{(t+u)^2}\,dt = \log u$$

Hence we obtain

$$\gamma_1(u) = \frac{1}{2}\left(\frac{1}{u}-1\right)\log u - \frac{1}{2}\log^2 u + \int_0^\infty \left[\psi'(t+u) - \frac{1}{t+u}\right]\log t\, dt \tag{3.13.2}$$

and with $u = 1$ we have

$$\gamma_1 = \int_0^\infty \left[\psi'(t+1) - \frac{1}{t+1}\right]\log t\, dt \tag{3.13.3}$$

and we see that (3.13.3) is in agreement with (2.4).

Differentiating (3.10) gives us

$$\int_0^\infty \frac{x\log\left(u^2+x^2\right)}{(u^2+x^2)(e^{2\pi x}-1)}\,dx - 2u\int_0^\infty \frac{\tan^{-1}\left(x/u\right)}{(u^2+x^2)(e^{2\pi x}-1)}\,dx$$

$$= -\int_0^\infty \left[\psi'(t+u) - \frac{1}{t+u} - \frac{1}{2(t+u)^2}\right]\log t\, dt$$

and substituting this in (3.13.1) results in

$$\gamma_1(u) = \frac{1}{2u}\log u - \frac{1}{2}\log^2 u + \int_0^\infty \frac{x\log\left(u^2+x^2\right)}{(u^2+x^2)(e^{2\pi x}-1)}\,dx - 2u\int_0^\infty \frac{\tan^{-1}\left(x/u\right)}{(u^2+x^2)(e^{2\pi x}-1)}\,dx \tag{3.13.3}$$

which concurs with the equivalent formula recently given for $\gamma_1$ by Choi [11].

□

We define $J(u)$ as

$$J(u) = \int_0^\infty \frac{\log\left(u^2+x^2\right)\tan^{-1}\left(x/u\right)}{e^{2\pi x}-1}\,dx$$

and the change of variable $x = y/2$ results in

$$J(u) = \frac{1}{2}\int_0^\infty \frac{\log\left(u^2+y^2/4\right)\tan^{-1}\left(y/2u\right)}{e^{\pi y}-1}\,dy$$

We see that

$$J(u/2)=\frac{1}{2}\int_0^\infty \frac{\left[\log\left(u^2+y^2\right)-\log 4\right]\tan^{-1}\left(y/u\right)}{e^{\pi y}-1}dy$$

$$=\frac{1}{2}\int_0^\infty \frac{\log\left(u^2+y^2\right)\tan^{-1}\left(y/u\right)}{e^{\pi y}-1}dy-\log 2\int_0^\infty \frac{\tan^{-1}\left(y/u\right)}{e^{\pi y}-1}dy$$

We define $K(u)$ as

$$K(u)=\int_0^\infty \frac{\tan^{-1}\left(x/u\right)}{e^{2\pi x}-1}dx$$

and, in view of (3.5), we have

$$K(u)=\frac{1}{2}\left[\log\Gamma(u)-\left(u-\frac{1}{2}\right)\log u+u-\frac{1}{2}\log(2\pi)\right]$$

The change of variable $x=y/2$ results in

$$K(u)=\frac{1}{2}\int_0^\infty \frac{\tan^{-1}\left(y/2u\right)}{e^{\pi y}-1}dy$$

and we easily see that

$$K(u/2)=\frac{1}{2}\int_0^\infty \frac{\tan^{-1}\left(y/u\right)}{e^{\pi y}-1}dy$$

We then have

$$\frac{1}{2}\int_0^\infty \frac{\log\left(u^2+y^2\right)\tan^{-1}\left(y/u\right)}{e^{\pi y}-1}dy=2K(u/2)\log 2+J(u/2)$$

Simple algebra shows us that

$$J(u)=\frac{1}{2}\int_0^\infty \frac{\log\left(u^2+x^2\right)\tan^{-1}\left(x/u\right)}{e^{\pi x}-1}dx-\frac{1}{2}\int_0^\infty \frac{\log\left(u^2+x^2\right)\tan^{-1}\left(x/u\right)}{e^{\pi x}+1}dx$$

$$=2K(u/2)\log 2+J(u/2)-\frac{1}{2}\int_0^\infty \frac{\log\left(u^2+x^2\right)\tan^{-1}\left(x/u\right)}{e^{\pi x}+1}dx$$

Therefore we have

$$H(u)=2[2K(u/2)\log 2+J(u/2)-J(u)]$$

where

$$H(u) \equiv \int_0^\infty \frac{\log\left(u^2+x^2\right)\tan^{-1}\left(x/u\right)}{e^{\pi x}+1}dx$$

From (3.10) we see that

$$J(u) = -\int_0^\infty \left[\psi(t+u) - \log(t+u) + \frac{1}{2(t+u)}\right]\log t\, dt$$

and we have the derivative

$$J'(u) = -\int_0^\infty \left[\psi'(t+u) - \frac{1}{t+u} - \frac{1}{2(t+u)^2}\right]\log t\, dt$$

Using (3.13.1) we see that

$$J'(u) = \frac{1}{2u}\log u - \frac{1}{2}\log^2 u - \gamma_1(u)$$

We then have

$$H'(u) = 2\left[K'(u/2)\log 2 + \frac{1}{2}J'(u/2) - J'(u)\right]$$

$$= 2\left\{K'(u/2)\log 2 + \frac{1}{2}\left[\frac{1}{u}\log(u/2) - \frac{1}{2}\log^2(u/2) - \gamma_1\left(\frac{u}{2}\right)\right] - \left[\frac{1}{2u}\log u - \frac{1}{2}\log^2 u - \gamma_1(u)\right]\right\}$$

Using $K'(u) = \frac{1}{2}\left[\psi(u) - \log u + \frac{1}{2u}\right]$ this becomes

$$H'(u) = \left[\psi\left(\frac{u}{2}\right) - \log(u/2) + \frac{1}{u}\right]\log 2 + \left[\frac{1}{u}\log(u/2) - \frac{1}{2}\log^2(u/2) - \gamma_1\left(\frac{u}{2}\right)\right]$$

$$-\left[\frac{1}{u}\log u - \log^2 u - 2\gamma_1(u)\right]$$

and we obtain

$$H'(u) = 2u\int_0^\infty \frac{\tan^{-1}\left(x/u\right)}{(u^2+x^2)(e^{\pi x}+1)}dx - \int_0^\infty \frac{x\log\left(u^2+x^2\right)}{(u^2+x^2)(e^{\pi x}+1)}dx$$

$$= \left[\psi\left(\frac{u}{2}\right) - \log(u/2) + \frac{1}{u}\right]\log 2 + \left[\frac{1}{u}\log(u/2) - \frac{1}{2}\log^2(u/2) - \gamma_1\left(\frac{u}{2}\right)\right]$$

$$-\left[\frac{1}{u}\log u-\log^2 u-2\gamma_1(u)\right]$$

In particular, we have

$$H'(1)=2\int_0^\infty \frac{\tan^{-1}x}{(1+x^2)(e^{\pi x}+1)}dx-\int_0^\infty \frac{x\log\left(1+x^2\right)}{(1+x^2)(e^{\pi x}+1)}dx$$

$$=\gamma_1+\gamma\log 2-\frac{1}{2}\log^2 2$$

where we have used (see for example[20])

$$\gamma_1\left(\frac{1}{2}\right)=\gamma_1-\log^2 2-2\gamma\log 2$$

One of the respondents to a question posed on the Mathematics Stack Exchange website came up with the following result in 2014

$$\int_0^\infty \frac{\tan^{-1}x}{(1+x^2)(e^{\pi x}+1)}dx=\frac{\pi^2}{16}-\frac{1}{4}-\frac{1}{4}\sum_{n=1}^\infty \frac{\log(n+1)}{n(n+1)}$$

and it may be noted that the latter series features prominently in [19].
http://math.stackexchange.com/questions/1056962/a-couple-of-definite-integrals-related-to-stieltjes-constants?rq=1

□

Making the substitution $x\to(t+u)x$ in (3.9) gives us

(3.14) $$\psi(t+u)-\log(t+u)+\frac{1}{2(t+u)}=-2\int_0^\infty \frac{x}{(1+x^2)(e^{2\pi(u+t)x}-1)}dx$$

where the parameter containing $(t+u)$ has thereby been switched from the quadratic in the denominator of the integral (3.9) to the exponential function in (3.14).

Integration with respect to $t$ gives us

$$I(u)=\int_0^\infty\left[\psi(t+u)-\log(t+u)+\frac{1}{2(t+u)}\right]dt=-2\int_0^\infty \frac{x}{1+x^2}\int_0^\infty \frac{1}{e^{2\pi(u+t)x}-1}dt$$

and we first of all consider the integral

$$J=\int_0^\infty \frac{1}{e^{2\pi(u+t)x}-1}dt=\int_0^\infty \frac{1}{e^{2\pi ux}e^{2\pi tx}-1}dt$$

Using $\frac{1}{y-1} = \frac{y}{y-1} - 1$ this becomes

$$J = \int_0^\infty \left[ \frac{e^{2\pi ux} e^{2\pi tx}}{e^{2\pi ux} e^{2\pi tx} - 1} - 1 \right] dt$$

$$= \frac{1}{2\pi x} \log(e^{2\pi ux} e^{2\pi tx} - 1) - t \Big|_0^\infty$$

We see that

$$\frac{1}{2\pi x} \log(e^{2\pi ux} e^{2\pi tx} - 1) - t = \frac{1}{2\pi x} \log[e^{2\pi tx} (e^{2\pi ux} - e^{-2\pi tx})] - t$$

$$= \frac{1}{2\pi x} \log(e^{2\pi ux} - e^{-2\pi tx})$$

and hence we have

$$J = u - \frac{1}{2\pi x} \log(e^{2\pi ux} - 1)$$

$$= -\frac{1}{2\pi x} \log(1 - e^{-2\pi ux})$$

We then have

$$I(u) = \frac{1}{\pi} \int_0^\infty \frac{\log(1 - e^{-2\pi ux})}{1 + x^2} dx$$

Integration by parts gives us

$$\int \frac{\log(1 - e^{-2\pi ux})}{1 + x^2} dx = \tan^{-1} x \log(1 - e^{-2\pi ux}) - \int \frac{2\pi u e^{-2\pi ux} \tan^{-1} x}{1 - e^{-2\pi ux}} dx$$

and hence we have the definite integral

$$(3.15) \qquad \int_0^\infty \frac{\log(1 - e^{-2\pi ux})}{1 + x^2} dx = -2u\pi \int_0^\infty \frac{\tan^{-1} x}{e^{2\pi ux} - 1} dx$$

Therefore we deduce that

$$I(u) = -2u \int_0^\infty \frac{\tan^{-1} x}{e^{2\pi ux} - 1} dx$$

$$= -2\int_0^\infty \frac{\tan^{-1}\left(t/u\right)}{e^{2\pi t}-1}\,dt$$

and we obtain

(3.16) $$I(u) = \int_0^\infty \left[\psi(t+u) - \log(t+u) + \frac{1}{2(t+u)}\right] dt = -2\int_0^\infty \frac{\tan^{-1}\left(x/u\right)}{e^{2\pi x}-1}\,dx$$

Hence substituting (3.5) we obtain

(3.17) $$\int_0^\infty \left[\psi(t+u) - \log(t+u) + \frac{1}{2(t+u)}\right] dt = \left(u - \frac{1}{2}\right)\log u - u + \frac{1}{2}\log(2\pi) - \log\Gamma(u)$$

which was also derived in [18].

Differentiating (3.17) gives us

$$\int_0^\infty \left[\psi'(t+u) - \frac{1}{t+u} - \frac{1}{2(t+u)^2}\right] dt = \log u - \frac{1}{2u} - \psi(u)$$

or equivalently

(3.18) $$\int_0^\infty \left[\psi'(t+u) - \frac{1}{t+u}\right] dt = \log u - \psi(u)$$

of which (2.2) is a particular case.

## 4. Another approach to the higher derivatives of the Riemann zeta function

Using Euler's reflection formula

$$\Gamma(s)\Gamma(1-s) = \frac{\pi}{\sin \pi s}$$

we may write the Riemann functional equation (1.21) as

$$s\varsigma(1-s) = 2(2\pi)^{-s}\,\frac{\Gamma(1+s/2)\Gamma(1-s/2)}{\Gamma(1-s)}\,\varsigma(s)$$

or equivalently as

$$2\varsigma(s) = \frac{1}{\lambda(s)}\,s\varsigma(1-s) \equiv \frac{1}{\lambda(s)}\,f(s)$$

where we have defined $\lambda(s)$ as

$$\lambda(s) = (2\pi)^{-s} \frac{\Gamma(1+s/2)\Gamma(1-s/2)}{\Gamma(1-s)}$$

and

$$f(s) = s\varsigma(1-s)$$

From the Laurent expansion of the Riemann zeta function

$$\varsigma(s) = \frac{1}{s-1} + \sum_{n=0}^{\infty} \frac{(-1)^n}{n!} \gamma_n (s-1)^n$$

we see that with $s \to 1-s$

$$s\varsigma(1-s) = -1 + \sum_{n=0}^{\infty} \frac{1}{n!} \gamma_n s^{n+1}$$

and hence we have

$$f^{(k)}(0) = k\gamma_{k-1}$$

We note that

$$\frac{d}{ds}\frac{1}{\lambda(s)} = -\frac{1}{\lambda(s)} g(s)$$

where

$$g(s) = \psi(1-s) + \frac{1}{2}\psi\left(1+\frac{s}{2}\right) - \frac{1}{2}\psi\left(1-\frac{s}{2}\right) - \log(2\pi)$$

with the result that the higher derivatives may be expressed as

$$\frac{d^i}{ds^i}\frac{1}{\lambda(s)} = \frac{1}{\lambda(s)} Y_i\left(-g(0), -g^{(1)}(0), ..., -g^{(i-1)}(0)\right)$$

where

$$g^{(i)}(s) = (-1)^i \psi^{(i)}(1-s) + \frac{1}{2^{i+1}}\psi^{(i)}\left(1+\frac{s}{2}\right) - \frac{1}{2^{i+1}}(-1)^i \psi^{(i)}\left(1-\frac{s}{2}\right)$$

We have

$$g^{(i)}(0) = \left[(-1)^i + \frac{1}{2^{i+1}}[1-(-1)^i]\right]\psi^{(i)}(1)$$

and substituting [28, p.22]

$$\psi^{(i)}(1) = (-1)^{i+1} i!\varsigma(i+1,1) = (-1)^{i+1} i!\varsigma(i+1)$$

we obtain

$$g^{(i)}(0) = \left[\frac{1}{2^{i+1}}[(-1)^{i+1}+1]-1\right] i!\varsigma(i+1)$$

where, for example, we have

$$g(0) = -[\gamma + \log(2\pi)]$$

$$g^{(1)}(0) = -\frac{1}{2}\varsigma(2)$$

Hence we have using the Leibniz differentiation rule

$$2\varsigma^{(n)}(0) = \sum_{i=0}^{n} \binom{n}{i} Y_i\left(-g(0), -g^{(1)}(0), ..., -g^{(i-1)}(0)\right) f^{(n-i)}(0)$$

and thus we obtain

(4.1) $$2\varsigma^{(n)}(0) = \sum_{i=0}^{n} \binom{n}{i} Y_i\left(-g(0), -g^{(1)}(0), ..., -g^{(i-1)}(0)\right)(n-i)\gamma_{n-i-1}$$

In passing we note that it was recently shown in [20] that this corresponds with

(4.2) $$n\gamma_{n-1} = 2\sum_{i=0}^{n} \binom{n}{i} Y_i\left(g(0), g^{(1)}(0), ..., g^{(i-1)}(0)\right)\varsigma^{(n-i)}(0)$$

The eta constants $\eta_n$ are defined by reference to the logarithmic derivative of the Riemann zeta function [13]

(4.3) $$\frac{d}{ds}[\log \varsigma(s)] = \frac{\varsigma'(s)}{\varsigma(s)} = -\frac{1}{s-1} - \sum_{k=0}^{\infty} \eta_k (s-1)^k \qquad |s-1| < 3$$

and, noting that $\lim_{s\to1}[(s-1)\varsigma(s)] = 1$, we obtain upon integration

(4.4) $$\log[(s-1)\varsigma(s)] = -\sum_{k=1}^{\infty} \frac{\eta_{k-1}}{k}(s-1)^k$$

We see from (4.3) that $\frac{\varsigma'(0)}{\varsigma(0)} = 1 - \sum_{k=0}^{\infty}(1)^k \eta_k$ and, since the series is convergent, we deduce that $\lim_{k\to\infty} \eta_k = 0$.

Coffey [13] has shown that the sequence $(\eta_n)$ has strict sign alteration, i.e.

(4.5) $$\eta_n = (-1)^{n+1}\varepsilon_n$$

where $\varepsilon_n$ are positive constants.

We showed in [16] that for $n \ge 0$

(4.6) $$(-1)^n (n+1)\gamma_n = Y_{n+1}\left(-0!\eta_0, -1!\eta_1, ..., -n!\eta_n\right)$$

and hence we obtain

$$2\varsigma^{(n)}(0) = \sum_{i=0}^{n} \binom{n}{i} Y_i\left(-g(0), -g^{(1)}(0), ..., -g^{(i-1)}(0)\right)(-1)^{n-i-1} Y_{n-i}\left(\gamma, -1!\eta_1, ..., -(n-i-1)!\eta_{n-i-1}\right)$$

We note that [10, p.412]

(4.7) $$Y_m(-x_1, x_2, ..., (-1)^m x_m) = (-1)^m Y_m(x_1, ..., x_m)$$

and hence the above equation may be expressed as

$$2\varsigma^{(n)}(0) = -\sum_{i=0}^{n} \binom{n}{i} Y_i\left(-g(0), -g^{(1)}(0), ..., -g^{(i-1)}(0)\right) Y_{n-i}\left(-\gamma, -1!\eta_1, ..., (-1)^{n-i}(n-i-1)!\eta_{n-i-1}\right)$$

We have [10, p.448]

(4.8) $$Y_n(x_1 + y_1, ..., x_n + y_n) = \sum_{i=0}^{n} \binom{n}{i} Y_i(y_1, ..., y_i) Y_{n-i}(x_1, ..., x_{n-i})$$

and this results in

$$2\varsigma^{(n)}(0) = -Y_n\left(-[g(0)+\gamma], -[g^{(1)}(0) + 1!\eta_1], ..., -[g^{(n-1)}(0) + (-1)^n (n-1)!\eta_{n-1}]\right)$$

which may be written as

$$2\varsigma^{(n)}(0) = -Y_n\left(\log(2\pi), \left[\frac{1}{2}\varsigma(2) - \eta_1\right]1!, ..., \left[\left\{1 - \frac{1}{2^n}[(-1)^n + 1]\right\}\varsigma(n) - (-1)^n \eta_{n-1}\right](n-1)!\right)$$

Using (4.7) again gives us

(4.9)
$$2\varsigma^{(n)}(0) = (-1)^{n+1} Y_n\left(-\log(2\pi), \left[\frac{1}{2}\varsigma(2) - \eta_1\right]1!, ..., \left[\left\{(-1)^n - \frac{1}{2^n}[(-1)^n + 1]\right\}\varsigma(n) - \eta_{n-1}\right](n-1)!\right)$$

and we designate $d_n$ as

$$d_n \equiv \left\{(-1)^n - \frac{1}{2^n}[(-1)^n + 1]\right\}\varsigma(n) - \eta_{n-1}$$

For example, with $n = 2$ in (4.9) we see that

$$2\varsigma''(0) = -Y_2\left(-\log(2\pi), \left[\frac{1}{2}\varsigma(2) - \eta_1\right]1!\right)$$

or equivalently

$$\varsigma''(0) = -\frac{1}{2}\log^2(2\pi) - \frac{1}{4}\varsigma(2) + \frac{1}{2}\eta_1$$

and since $\eta_1 = \gamma^2 + 2\gamma_1$ this concurs with (1.10).

The $\sigma_k$ constants are defined by the Taylor expansion

(4.10) $$\log \xi(s) = -\log 2 - \sum_{k=1}^{\infty} \frac{\sigma_k}{k} s^k$$

where $\xi(s)$ is the Riemann xi function defined by

$$\xi(s) = \frac{1}{2} s(s-1)\pi^{-s/2}\Gamma(s/2)\varsigma(s)$$

It was shown by Zhang and Williams in 1994 that (see for example [16] and the references therein)

(4.11) $$\eta_n = (-1)^{n+1}\left[\sigma_{n+1} + \left(1 - \frac{1}{2^{n+1}}\right)\varsigma(n+1) - 1\right]$$

and hence we have

(4.12) $$d_n = (-1)^n - (-1)^n\left[\sigma_n + \frac{(-1)^n}{2^n}\varsigma(n)\right]$$

Lehmer [27] considered the constants $b_n$ defined by

(4.13) $$\frac{d}{ds}\log[2(s-1)\varsigma(s)] = \frac{\varsigma'(s)}{\varsigma(s)} + \frac{1}{s-1} = \sum_{n=0}^{\infty} b_n s^n \quad , \ |s| < 2$$

so that

(4.14) $$\log[2(s-1)\varsigma(s)] = \sum_{n=0}^{\infty} \frac{b_n}{n+1} s^{n+1} = \sum_{n=1}^{\infty} \frac{b_{n-1}}{n} s^n$$

We note that

(4.15) $$\frac{d^m}{ds^m}[2(s-1)\varsigma(s)]\bigg|_{s=0} = Y_m(0!b_0, 1!b_1, ..., (m-1)!b_{m-1})$$

(4.16) $$2\left[m\varsigma^{(m-1)}(0) - \varsigma^{(m)}(0)\right] = Y_m(0!b_0, 1!b_1, ..., (m-1)!b_{m-1})$$

Lehmer [27] showed that

(4.17) $$b_{n-1} = -\left[\sigma_n + \frac{(-1)^n}{2^n}\varsigma(n)\right]$$

and accordingly we obtain from (4.12)

(4.18) $$d_n = (-1)^n[1 + b_{n-1}]$$

Hence we have from (4.9)

$$2\varsigma^{(n)}(0) = (-1)^{n+1}Y_n\left(-\log(2\pi), [1+b_1]1!, ..., (-1)^n[1+b_{n-1}](n-1)!\right)$$

Since [27] $b_0 = \log(2\pi) - 1$ we have

(4.19) $$2\varsigma^{(n)}(0) = (-1)^{n+1}Y_n\left(-[1+b_0]0!, [1+b_1]1!, ..., (-1)^n[1+b_{n-1}](n-1)!\right)$$

Coffey [13] showed in 2008 that $b_m$ has strict sign alteration, i.e.

(4.20) $$b_m = (-1)^m \mu_m \text{ where } \mu_m > 0$$

This strict sign alteration was also independently reported in [16] in 2009 where we considered the function $L(s) = \log[(s-1)\varsigma(s)]$. Hence we have

$$2\varsigma^{(n)}(0) = (-1)^{n+1}Y_n\left(-[1+\mu_0]0!, [1-\mu_1]1!, ..., (-1)^n[1+(-1)^{n-1}\mu_{n-1}](n-1)!\right)$$

and employing (4.7) this may be expressed as

(4.21) $$2\varsigma^{(n)}(0) = -Y_n\left([1+\mu_0]0!, [1-\mu_1]1!, ..., [1+(-1)^{n-1}\mu_{n-1}](n-1)!\right)$$

It is clear that a sufficient (but not necessary) condition for $\varsigma^{(n)}(0)$ to be negative is that $\mu_n < 1$ for all $n \geq 0$ because all of the arguments of the complete Bell polynomial $Y_n\left([1+\mu_0]0!, [1-\mu_1]1!, ..., [1+(-1)^{n-1}\mu_{n-1}](n-1)!\right)$ would then be positive. We do know however that $\lim_{n\to\infty} b_n = 0$ and hence there exists an $N$ such that $|b_n| < 1$ for all $n \geq N$. This fact that $\varsigma^{(n)}(0)$ is negative for sufficiently large values of $n$ may also be deduced more readily from the known limit $\lim_{n\to\infty}\frac{\varsigma^{(n)}(0)}{n!} = -1$ noted by Apostol [4].

We recall (4.11)

$$\eta_n = (-1)^{n+1}\left[\sigma_{n+1} + \left(1 - \frac{1}{2^{n+1}}\right)\varsigma(n+1) - 1\right]$$

and using (4.5) $\eta_n = (-1)^{n+1}\varepsilon_n$ we see that

$$\varepsilon_n = \sigma_{n+1} + \left(1 - \frac{1}{2^{n+1}}\right)\varsigma(n+1) - 1$$

and hence we have

(4.22) $$\sigma_{n+1} > 1 - \left(1 - \frac{1}{2^{n+1}}\right)\varsigma(n+1)$$

We recall Lehmer's relation (4.17) for $n \geq 1$

$$\sigma_{n+1} = (-1)^n 2^{-n-1}\varsigma(n+1) - b_n$$

and we deduce that for $n \geq 1$

(4.23) $$\varsigma(n+1) - 1 - [1 + (-1)^{n+1}]\frac{\varsigma(n+1)}{2^{n+1}} > b_n$$

With $n \to 2n$ this inequality becomes

$$\varsigma(2n+1) - 1 > b_{2n}$$

Since $\varsigma'(s) < 0$ for all $s > 1$ we deduce that $\varsigma(s)$ is monotonic decreasing for all $s > 1$. We therefore have $\varsigma(3) > \varsigma(2n+1)$ for all $n \geq 1$ and, since $\varsigma(3) \simeq 1.2020...$, it is easily seen that $1 > \varsigma(3) - 1 > \varsigma(2n+1) - 1$ and hence for all $n \geq 1$ we have

$$1 > b_{2n}$$

With $n \to 2n-1$ in (4.23) we obtain

$$\varsigma(2n)\left[1 - \frac{1}{2^{2n}}\right] - 1 > b_{2n-1}$$

but unfortunately this does not appear to assist us in trying to prove that $1 > |b_{2n-1}|$.

□

We may obtain an expression for $b_n$ in terms of $\varsigma^{(n)}(0)$ by reference to the following inversion relation of Chou et al. [12]

$$Y_n = y_n = \sum_{k=1}^{n} B_{n,k}\left(x_1, x_2, ..., x_{n-k+1}\right) \Leftrightarrow x_n = \sum_{k=1}^{n} (-1)^{k-1}(k-1)!B_{n,k}\left(y_1, y_2, ..., y_{n-k+1}\right)$$

We note that [10, p.415]

$$Y_{n+1}(x_1, ..., x_n) = \sum_{i=0}^{n} \binom{n}{i} Y_{n-i}(x_1, ..., x_{n-i}) x_{i+1}$$

and using (4.19) we obtain the recurrence

(4.24) $$\varsigma^{(n+1)}(0) = \sum_{i=0}^{n} \binom{n}{i} (i)!\varsigma^{(n-i)}(0)[1+b_i]$$

**5. Open access to our own work**

This paper contains references to various other papers and, rather surprisingly, most of them are currently freely available on the internet. Surely now is the time that all of our work should be freely accessible by all. The mathematics community should lead the way on this by publishing everything on arXiv, or in an equivalent open access repository. We think it, we write it, so why hide it? You know it makes sense.

Wessex House,
Devizes Road,
Upavon,
Pewsey,
Wiltshire SN9 6DL
dconnon@btopenworld.com